\documentclass[11pt]{amsart}
\usepackage{amsmath}
\usepackage{amscd}
\usepackage{amsfonts,latexsym,amssymb}

\renewcommand{\b}{\beta}

\renewcommand{\l}{\lambda}
\newcommand{\s}{\sigma}

\renewcommand{\o}{\omega}
\newcommand{\p}{\phi}

\newcommand{\ima}{\mathbf{i}}

\newcommand{\G}{\Gamma}
\renewcommand{\O}{\Omega}

\newcommand{\SC}{{\mathcal{C}}}

\newcommand{\SF}{{\mathcal{F}}}

\newcommand{\SL}{{\mathcal{L}}}
\newcommand{\SM}{{\mathcal{M}}}

\newcommand{\SO}{{\mathcal{O}}}

\newcommand{\SU}{{\mathcal{U}}}

\newcommand{\Z}{\mathbb{Z}}
\newcommand{\C}{\mathbb{C}}
\newcommand{\Q}{\mathbb{Q}}

\newcommand{\SSS}{\mathbb{S}}

\newcommand{\isom}{\cong}

\newcommand{\Jac}{\operatorname{Jac}}

\newcommand{\Gr}{\operatorname{Gr}}
\newcommand{\rank}{\operatorname{rank}}

\newcommand{\Hom}{\operatorname{Hom}}

\newcommand{\dimC}{\text{dim}_{\C}\,}

\newcommand{\inc}{\hookrightarrow}

\newcommand{\End}{\operatorname{End}}
\newcommand{\Tr}{\operatorname{Tr}}

\newcommand{\la}{\langle}
\newcommand{\ra}{\rangle}
\newcommand{\bd}{\partial}
\newcommand{\bbd}{\bar{\partial}}
\newcommand{\x}{\times}
\newcommand{\ox}{\otimes}

\newtheorem{proposition}{Proposition}[section]
\newtheorem{theorem}[proposition]{Theorem}
\newtheorem{definition}[proposition]{Definition}

\newtheorem{corollary}[proposition]{Corollary}

\title{Semipositive bundles and Brill-Noether theory}

\thanks{Partially supported by The European Contract Human
Potential Programme, Research Training Network
HPRN-CT-2000-00101.}

\date{July, 2001.}


\author{Vicente Mu\~noz}
\address{Departamento de Matem\'aticas \\
Facultad de Ciencias \\ Universidad Aut\'onoma de Madrid
\\ 28049 Madrid \\ Spain}
 \email{vicente.munoz@uam.es}

\author{Francisco Presas}
\address{Departamento de \'Algebra \\ Facultad de
Ciencias Matem\'aticas \\ Universidad Complutense de Madrid \\
28040 Madrid \\ Spain} \email{fpm@eucmos.sim.ucm.es}

\begin{document}

\begin{abstract}
We prove a Lefschetz hyperplane theorem for the determinantal loci
of a morphism between two holomorphic vector bundles $E$ and $F$
over a complex manifold under the condition that $E^*\ox F$ is
Griffiths $k$-positive. We apply this result to find some homotopy
groups of the Brill-Noether loci for a generic curve.
\end{abstract}

\maketitle

\section{Introduction}
\label{introduction}

The topological properties of the zero sets of sections of a very
ample line bundle on a projective variety were studied by
Lefschetz in the 20's, when he proved the renowned Lefschetz
hyperplane theorem. Later on Bott, Andreotti and Fraenkel
\cite{Bott, AF} gave a new proof of the result using Morse theory,
which has the byproduct of giving homotopy isomorphisms. The
ampleness condition of the vectors bundles in these proofs can be
weakened to a {\em semipositivity\/} condition. We recall the
concept of Griffiths $k$-positivity \cite{nakano}

\begin{definition} \label{k_positive}
A holomorphic vector bundle $E$ over a complex manifold $M$ is
said to be \textup{(}Griffiths\textup{)} $k$-positive if there
exists a hermitian metric $h$ on $E$ such that for every point
$x\in M$ there exists a complex subspace $V_x\subset T_xM$ of
dimension at least $k$ where the curvature form $\Theta$ of the
connection associated to $h$ satisfies that $\Theta_{v,iv}$ is a
definite positive quadratic form in the fiber $E_x$, for every
non-zero $v\in V_x$.
\end{definition}

For $k=\dim_{\C} M$ we recover the definition of positivity of a
vector bundle \cite{Gr69}. We have the following extension of the
classical Lefschetz theorem

\begin{theorem} \label{main_thm}
Let $E$ be a rank $r$ $k$-positive vector bundle over a complex
compact manifold $M$ and let $s$ be a holomorphic section of $E$.
Let $W=Z(s)$ be the zero set of $s$. Then $M-W$ has the homotopy
type of a CW-complex of dimension $2n-(k-r)-1$.
When $W$ is smooth the inclusion $W \to M$ induces isomorphisms on
homology \textup{(}resp.\ homotopy\textup{)} groups $H_p$
\textup{(}resp.\ $\pi_p$\textup{)} for $p<k-r$ and an epimorphism
for $p=k-r$.
\end{theorem}

Theorem \ref{main_thm} may be extended in the following way.

\begin{definition}\label{def:kpositivepair}
Let $E$ be a holomorphic vector bundle over a complex manifold $M$
and let $s$ be a holomorphic section of $E$. Then the pair $(E,s)$
is $k$-positive if there exists a hermitian metric $h$ on $E$ such
that for every $x\in M$ which is a critical point for $|s|^2$ with
$s(x)\neq 0$, there is a complex subspace $V_x\subset T_xM$ of
dimension at least $k$ satisfying $h(s, \Theta_{v,iv}s)
>0$ for every non-zero $v\in V_x$.
\end{definition}

\begin{theorem}\label{thm:main2}
 Let $E$ be a rank $r$ vector bundle over a compact
 complex manifold $M$ and $s$ a holomorphic section such that
 $(E,s)$ is a $k$-positive pair. Let $W=Z(s)$ be the zero set of
 $s$. Then $M-W$ has the homotopy type of a CW-complex of dimension
 $2n-(k-r)-1$.
 When $W$ is smooth the inclusion $W \to M$ induces isomorphisms on
 homology \textup{(}resp.\ homotopy\textup{)} groups $H_p$
 \textup{(}resp.\ $\pi_p$\textup{)} for
 $p<k-r$ and an epimorphism for $p=k-r$.
\end{theorem}

Let us suppose now that we have two holomorphic vector bundles $E$
and $F$ over $M$ of ranks $e$ and $f$ respectively, and let
$\phi:E\to F$ be a vector bundle morphism. We can always assume
that $e\leq f$ by changing to the transpose morphism if necessary.
Fix a positive integer $r\leq e$. The $r$-determinantal subvariety
of $\phi$ is defined as
 $$
  D_r(\phi)= \{ x\in M: \rank \, \phi_x \leq r\}.
 $$
This can be constructed in the following way. Consider the
grassmannian fibration $\pi: G=\Gr(e-r, E)\to M$, and the
composition
 $$
  \phi_r:U \hookrightarrow \pi^* E
  \stackrel{\pi^*\phi}{\longrightarrow}\pi^* F,
 $$
where the first map is the natural inclusion of the universal
bundle $U$ over $G$. The zero set of $\phi_r$ satisfies
 $$
  \pi(Z(\phi_r))=D_r(\phi).
 $$
In fact $\pi$ is bijective over $D_r(\phi)-D_{r-1}(\phi)$. We have
the following result which extends that of \cite{FL, Debarre}

\begin{theorem}\label{thm:determ}
Let $\phi$ be a morphism between the holomorphic vector bundles
$E$ and $F$ over a compact complex manifold $M$, such that
$E^*\otimes F$ is $k$-positive. Suppose that $Z(\p_r)$ is smooth.
Then the natural inclusion $Z(\phi_r) \hookrightarrow
G=\Gr(e-r,E)$ induces isomorphisms between the homology
\textup{(}resp.\ homotopy\textup{)} groups $H_i(Z(\phi_r))$ and
$H_i(G)$ \textup{(}resp.\ $\pi_i(Z(\phi_r))$ and
$\pi_i(G)$\textup{)} for $i< k-(e-r)(f-r)$ and epimorphisms for
$i= k-(e-r)(f-r)$.
\end{theorem}

We apply the theory developed to recover Debarre's result
\cite{Debarre} computing the homology of Brill-Noether loci over
an algebraic curve. Our method gives us as well information about
the homotopy groups. Let $C$ be an algebraic curve. Consider the
jacobian variety of degree $d$ line bundles $\Jac^d(C)$. The
Brill-Noether moduli spaces, $W_d^k(C)$, are defined as
 $$
W_d^k(C)= \{ L\in \Jac^d(C)| ~ h^0(L)\geq k+1 \}.
 $$
The varieties of linear systems are
 $$
G_d^k(C)= \{ (L,V)| L \in W_d^k(C), V\subset H^0(L), \dim V = k+1
\}.
 $$
Clearly there is a map $G_d^k(C) \to W_d^k(C)$ which is one-to-one
over $W_d^k(C)-W^{k+1}_d(C)$. For a generic curve $C$, the spaces
$G_d^k(C)$ are smooth \cite{ACGH85} and of the expected dimension
$\rho=g-(k+1)(g-d+k)$. We prove the following results about the
topology of $G_d^k(C)$ and $W_d^k(C)$.

\begin{theorem} \label{tue}
Let $C$ be a generic curve of genus $g\geq 2$ and let $d\geq 1$
and $k\geq 0$ be integers. If $\rho>0$ then $G_d^k(C)$ and
$W_d^r(C)$ are connected and non-empty. If $\rho>1$ then
$\pi_1(W_d^r(C))=\pi_1(G_d^r(C))$ is free abelian of rank $2g$. In
general we have that
 $$
  \pi_i(G_d^k(C))=\pi_i(BU(k+1)),\, 2\leq i\leq \rho-1.
 $$
In particular, the corresponding rational homotopy groups of
$G_d^k(C)$ are
 $$
  \pi_i(G_d^k(C)) \ox \Q= \left\{ \begin{array}{ll} \Q \qquad &
  i \leq 2k+2, \> i\equiv 0 \pmod2 \\  0 &\text{otherwise}
  \end{array} \right.
 $$
 for $2\leq i\leq \rho-1$.
\end{theorem}

\begin{corollary} \label{tue2}
  Let $C$ be a generic curve of genus $g\geq 2$ and let $d\geq 1$,
  $k\geq 0$ and $l\geq 1$ be integers.
  Suppose that $G_d^{k+i}(C)$ are non-empty for $0\leq i\leq l$. Then
  $\pi_1(W_d^k(C))=\Z^{2g}$ and $\pi_i(W_d^k(C))=0$ for $2\leq i\leq 2l$.
\end{corollary}

{\em Acknowledgements:} We are specially indebted to Ignacio Sols
for proposing us this problem. We are also grateful to Ignasi
Mundet and Aniceto Murillo for very useful comments.

\section{Proof of Theorems \ref{main_thm} and \ref{thm:main2}}
\label{proof}

The result of theorem \ref{main_thm} follows from that of theorem
\ref{thm:main2} since  if $E$ is $k$-positive then $(E,s)$ is a
$k$-positive pair for any holomorphic section $s$. Therefore we
may suppose that we are under the assumptions of theorem
\ref{thm:main2}. Then the result is a slight generalization of the
classical proof using Morse theory. We consider the real function
$f= \log h(s,s)$. Clearly $f$ is well defined over $M-Z(s)$ where
$Z(s)$ is the zero set of $s$. Our objective is to show that all
the critical points of $f$ have index at least $k-r+1$, since then
the result follows from standard Morse theory.

Now we compute
 $$
  \bd f= \frac{1}{h(s,s)} (h(\bbd s, s) + h(s, \bd s)).
 $$
At a critical point this quantity vanishes and so we get
 $$
  h( \bbd s, s) + h(s, \bd s)=0.
 $$
Recalling that $\bbd s=0$ we obtain
\begin{equation}
 h(s,\bd s) = 0. \label{null}
\end{equation}
A second differentiation, omitting quantities that vanish at the
critical point, gives us
 $$
 \bbd \bd \log |s|^2= \frac{1}{|s|^2} (h(\bd \bbd s, s) +
 h( \bbd s, \bbd s) + h( \bd s, \bd s) + h(s, \bbd \bd s)).
 $$
Now observe that the section is holomorphic and so $\bbd s=0$.
Also we have that $\bbd \bd + \bd \bbd = \Theta$, the curvature of
the bundle $E$. Substituting into the equation we obtain
 $$
 \bbd \bd \log |s|^2= \frac{1}{|s|^2} (h(s,\Theta s) + h( \bd s, \bd s)).
 $$
Now we must control the second term in the sum. For this we define
the subspace
 $$
 W= \{ v\in V_x | ~\bd s(v)=0 \}.
 $$
Using that for any $v\in T_xM$ we have that $h(\bd s(v),s)=0$, we
have that the complex dimension of $W$ is bounded below by
$k-r+1$. So $h(\bd s, \bd s)=0$  in $W$. By the condition on the
curvature, it follows that, restricting to $v\in W$, $\bbd \bd
f(v, \ima v) = -\ima \alpha$, for some real number $\alpha>0$. Let
$H_f$ be the Hessian of $f$. Then
 $$
  H_f(v) + H_f(\ima v) = -2\ima \bbd \bd f(v,  \ima v).
 $$
This quantity is clearly negative for any non-zero $v\in W$.
Suppose that the index of the critical point were strictly less
than $k-r+1$. Then there would be a subspace $P\subset T_xM$ of
dimension greater than $2n-k+r-1$, where $H_f$ is definite
positive. Here $n=\dimC M$. So $P \bigcap \ima P$ would have real
dimension at least $2n-2k+2r$. This subspace would intersect
non-trivially to $W$ providing a contradiction. So the index of
the critical point is at least $k- r+1$.

Now a standard argument in Morse theory gives that $M-W$ has the
homotopy type of a CW-complex of dimension $2n-(k-r+1)$.

When $W$ is smooth, we consider a small tubular neighbourhood
$E(W)$ that retracts to $W$. Then $M$ is obtained from $E(W)$ by
attaching cells of index at least $k-r+1$. This means that the
$(k-r)$-skeleta of $M$ and $W$ are the same (up to homotopy).
Therefore we get that the natural map $\pi_p(W)\to \pi_p(M)$ is an
isomorphism for $p<k-r$, and an epimorphism for $p=k-r$. The
statement about homology groups holds by the same reason.
\hfill $\Box$ \vspace{5mm}

\section{Determinantal submanifolds} \label{determ}

We are going to apply our results in the previous section to the
study of determinantal subvarieties, aiming to prove theorem
\ref{thm:determ}. Suppose we have a morphism $\phi$ between
bundles $E$ and $F$ over a complex manifold $M$. We suppose that
$\rank\, E=e\leq f=\rank\, F$. The morphism $\phi$ can be
interpreted as a holomorphic section of the bundle $E^*\ox F$. Fix
a positive integer $r\leq e$, we can define the manifold
$G=\Gr(e-r, E)$ which is a grassmannian fibration over $M$. The
canonical projection of this fibration will be denoted by $\pi$.
We consider the morphism $\phi_r:U \inc \pi^* E \to \pi^* F$ on
$G$, where $U$ is the universal bundle. The $r$-determinantal
subvariety of $\phi$ is
 $$
 D_r(\phi)=\pi(Z(\phi_r))= \{ x\in M: \rank \, \phi_x \leq r \}.
 $$
The subvarieties $\{ D_r(\phi)\}_{i=0}^{e-1}$ have a natural
structure of a stratified submanifold. In some situations, the
strata are irreducible and reduced complex varieties of the
expected dimension. In this case, $Z(\phi_r)$ are smooth
subvarieties of $G$, so that they can be considered as a
desingularization of $D_r(\phi)$. Note that $\dim G=\dim M+
r(e-r)$ and $\rank (U^*\ox \pi^*F)= f(e-r)$, so that the expected
(complex) codimension of $Z(\phi_r)$ is $(e-r)(f-r)$. If this
exceeds the dimension of $M$ then $D_r(\phi)$ is empty.

We start by relating the curvature of the bundles over $M$ and
over $G$.

\begin{proposition}\label{G}
Let $E$ be a rank $e$ $k$-positive vector bundle over $M$. The
vector bundle $U$ on $G=\Gr(e-r,E)$ is $(k+r)$-positive.
\end{proposition}

{\bf Proof.} Let $x\in M$ and choose holomorphic coordinates
$(z_1, \ldots, z_n)$ at a neighborhood of the point. Fix a point
$V\in G$ such that $\pi(V)=x$. Denote $p=e-r$. Choose a
holomorphic frame $f=(u_1, \ldots u_e)$ of the bundle $E$ at a
neighborhood of $x$. We can assume that $u_1, \ldots u_p$ expand
the subspace $V$ and also that $f(0)$ is an orthonormal basis with
respect to the hermitian metric $h$ of $E$.

With these choices it is now easy to define holomorphic
coordinates at a neighborhood of $V\in G$. The chart is defined as
\begin{eqnarray*}
  \C^n \times \C^{p(e-p)} & \to & G \\
  (z_l,\l_{kj}) & \to & < u_1 + \sum_{ k=p+1}^e \l_{k1}
  u_k, \ldots, u_p + \sum_{k=p+1}^e \l_{kp} u_k>.
\end{eqnarray*}
A chart for the bundle $U$ is given by
\begin{eqnarray} \label{chart}
 \C^n \times \C^{p(e-p)} \times \C^p & \to & U\subset \pi^* E \\
 (z_l, \l_{kj}, \beta_j) & \to & (u_1 + \sum_{k=p+1}^e \l_{k1} u_k,
 \ldots, u_p + \sum_{k=p+1}^e \l_{kp} u_k) \left(
 \begin{array}{c} \beta_1 \\ \vdots \\ \beta_p
 \end{array} \right). \nonumber
\end{eqnarray}
With this chart the pull-back of the metric $h$ to $\pi^* E$ is
defined by the matrix $h(f)=(h(u_i, u_j))_{i,j}$. This matrix can
be divided in the following blocks
 $$ h(f)= \left(
\begin{array}{cc} h^{11}(f) & h^{12}(f) \\ h^{21}(f) & h^{22}(f)
\end{array} \right),
 $$
according to the local decomposition $E= <u_1, \ldots, u_p>
\bigoplus <u_{p+1}, \ldots, u_e >$.

For a point of $G$ with coordinates $(z_l, \l_{kj})$ the
restriction of the metric $h$ to $U$ with respect to the frame
$f_U=(u_1 + \sum_{k=p+1}^e \l_{k1} u_k, \ldots, u_p +
\sum_{k=p+1}^e \l_{kp} u_k)$ is given by
 $$
 h_U(f_U)= h^{11}(f) +
h^{12}(f) \l + \l^* h^{21}(f) + \l^* h^{22}(f) \l,
 $$
where $\l^*$ is the conjugate transpose matrix of $\l$. If we use
a holomorphic frame $f$ such that $h(z)= I + O(|z|^2)$ then we can
compute the curvature using Lemma 2.3 in Chap. III of \cite{We73}.
So the curvature in this case is $\Theta(0)= \bbd\bd h(0)$. We
start by computing the derivative of $h_U$ (we do not write the
frame $f$ to simplify the notation),
 $$
  \bd h= \bd_z h^{11} + \bd_z h^{12} \l +
  h^{12}d\l+ \l^* \bd_z h^{21} + \l^* \bd_z h^{22} \l+
  \l^* h^{22}d\l,
 $$
where by $\bd_z$ we denote the derivatives in the $(z_l)$
directions. Differentiating again
\begin{eqnarray*}
 \bbd \bd h & = & \bbd_z \bd_z h^{11} +
 \bbd_z \bd_z h^{12} \l + \bbd_z h^{12} \wedge d\l+
 d\l^* \wedge \bd_z h^{21} +\\
 & & + \l^* \bbd_z \bd_z h^{21} +  d\l^* \wedge
 \bd_z h^{22} \l  + \l^* \bbd_z \bd_z
 h^{22} \l + \l^* \bbd_z h^{22} \wedge d\l + d\l^*\wedge
 h^{22} d\l.
\end{eqnarray*}
At $\l=0$ we have $h^{11}=I_p$, $h^{22}=I_{r-p}$, $h^{12}=0$,
$h^{21}=0$ and $\bd h= \bbd h=0$. Therefore
\begin{equation}
 \Theta_U(f_U(0))= \bbd \bd h^{11}(f)(0) + d \l^* \wedge d \l=
 \Theta(0)_{|U} + d \l^* \wedge d \l. \label{curv_up}
\end{equation}
Now we study the grassmannian direction. We want to find a
subspace in $<\frac{\bd}{\bd \l_{kj}}>$ where $d \l^* \wedge d\l$
is positive definite for any non trivial pair of directions $v,iv$
in the base. So we need to impose that for any non-zero $\beta\in
\C^p$,
 \begin{equation}
   \b^* d \l^* \wedge d\l\, \b>0. \label{conditio}
 \end{equation}
We may choose $\Delta$ to be the subspace of those $\l\in
\C^{p(e-p)}$ which have all their columns
$\l_j=(\l_{kj})_{k=p+1}^e$ equal. For any $\b\neq 0$ and any
non-zero $\l \in \Delta$ we have
 $$
 \b^* d \l^* \wedge d\l\, \b= \sum_j |\b_j|^2 d\l_j^*\wedge d\l_j>0.
 $$

Take the subspace $H\subset T_xM$ where $\Theta$ is positive
definite at $x$. Then the subspace $H \bigoplus \Delta$ is a
subspace where $\Theta_U$ is positive definite. Recall that
$r=e-p$. So we have proved that $U$ is $(k+r)$-positive. It is
easy to check that we have found the largest dimension for a
subspace $\Delta$ satisfying the required properties. \hfill
$\Box$ \vspace{5mm}

Using proposition \ref{G} together with theorem \ref{main_thm}
would not give a large range for Lefschetz isomorphisms for
$Z(\phi_r)\subset G$. We get around this problem by using the
notion of $k$-positive pair in definition \ref{def:kpositivepair}
and the freedom of varying the hermitian metrics in $E$ and $F$.

\begin{proposition} \label{G2}
  Let $\p:E\to F$ be a morphism where $E^*\ox F$ is k-positive and
  let $\p_r:U\to \pi^*F$ be the induced map on $G=\Gr(e-r,E)$.
  Then $(U^*\ox \pi^* F, \p_r)$ is a $(k+r(e-r))$-positive pair.
\end{proposition}

{\bf Proof.}
  Let $V\in G$ be a critical point of $|\p_r|^2$. As $\p_r$ is
  holomorphic, this implies that
  $h(\p_r, \nabla \p_r)=0$ at $V$, where $h$ is the metric
  induced in $U^*\ox F$ by the metrics in $E$ and $F$.
  Let $x=\pi (V)\in M$ and choose
  holomorphic coordinates $z_l$ in $M$ with $z_1=\cdots=z_n=0$ at
  $x$. Fix holomorphic frames $f_E=(u_1,\ldots, u_e)$ and
  $f_F=(v_1,\ldots, v_f)$ for $E$ and $F$ respectively,
  orthonormal at $x$, such that $V=<u_1,\ldots, u_p>$,
  where $p=e-r$. The map $\p:E\to F$
  becomes, with respect to these frames, a map $s:\C^e\to \C^f$,
  which is decomposed as $s=[{s_1\atop s_2}]$,
  with respect to $\C^e=\C^p\oplus \C^r$.

  As in the proof of Proposition \ref{G} we have coordinates
  $(z_l,\l_{kj})$ for $G$ using the frame $f_E$. The bundle $U$
  has a natural trivialization given by \eqref{chart} and $\p_r:\C^p\to \C^f$
  is written as the map $\p_r=s_1 +s_2 \l$. Therefore $\nabla
  \p_r=\nabla_z s_1 + s_2 d\l$, at the point $V\in G$ (since this
  point has $\l=0$).
  The condition $h(\p_r, \nabla \p_r)=0$ translates thus into
  $$
    h(s_1, \nabla_z s_1)=0 \qquad \text{and} \qquad s_1^*s_2=0.
  $$

  We need to compute $h(\s,\Theta_{v,iv} \s)$, where $\Theta$ is the
  curvature form of $U^*\otimes F$ and the section is $\s=\p_r$.
  By \eqref{curv_up} we have
  $$
   \Theta(0)= \Theta_{E^*\ox F}(0) |_{U^*\ox F} -  d\l^* \wedge
   d\l.
  $$
  Let $v=(u,W)\in T_VG=T_xM\oplus T_V\Gr(p,e)$ be a vector. We
  have
  $$
   h(s_1, \Theta(0)_{v,iv}s_1)= h(s_1,\Theta_{E^*\ox
   F}(0)_{u,iu}s_1) - \ima \Tr (s_1^* W^* W s_1),
  $$
  by the usual formula for the metric in the tensor product $E^*
  \ox F$.
  If we can arrange that $s_1$ is injective then
  we would have $W s_1\neq 0$ for any $W\neq 0$. Then we choose $H\subset
  T_xM$ where $\Theta_{E^*\ox F}(0)_{u,iu}$ is definite positive and
  hence $\Theta_{v,iv}$ is positive on $H \oplus T_{V} \Gr(p,e)$,
  which has dimension $k+ (e-p)p=k+r(e-r)$.

  So it remains to prove that for suitable hermitian metrics on $E$
  and $F$ we have that $h(s_1,\nabla s_1)=0$ and $s_1^*s_2=0$
  imply that $s_1$ is injective. For this we use a perturbation
  argument. Fix the metric on $E$ and let
  $\SM_F$ be the space of hermitian metrics for $F$, completed in a
  suitable Sobolev norm so that it is a Hilbert manifold and the
  elements $g\in \SM_F$ have enough regularity. We define the map
  \begin{eqnarray*}
    \SF=(\SF_1,\SF_2): G \x \SM_F & \to& T^*M \oplus \Hom (U, U^{\perp}) \\
    (y, g) & \mapsto & \left( g(s_1,\nabla s_1)_{\pi(y)}, (s_2^{*_g}s_1)_y \right),
  \end{eqnarray*}
  where the dual of $s_2$ is defined using $g$.
  We want to check that this map is a submersion. Fix $(y,g_0)$ and
  let us compute the differential $D\SF$ at this point.
  Fix a trivialization of $G$
  as in the proof of proposition \ref{G} with coordinates
  $(z,\l)$. Fix a hermitian trivialization of $E$ (recall that we
  are only varying the metric of $F$) and a holomorphic
  trivialization of $F$ such that $g_0(0)=I$.
  The points in $G\x \SM_F$ are denoted by
  $(z,\l, g)$, where $g$ is a hermitian matrix valued function of
  $z$. The tangent vectors will be denoted as $(v,A,b)$, where $b$ is also a
  hermitian matrix valued function. The first component of $\SF$
  is $\SF_1=\bd |s_1|_g^2$, which in the chart gives
  $\SF_1(0,0,g)=\Tr (s_1^* g ds_1 + s_1^*\bd g s_1)$. Thus
  $$
  D\SF_1(0,0,b)=\Tr(s_1^* b(0) ds_1 + s_1^* (\bd b)(0) s_1).
  $$
  Also $\SF_2(0,\l,h)=(s_2-s_1\l^*)^* g (s_1+s_2\l)$ and then
  $$
  D\SF_2(0,A,b)=\Tr(s_2^*b(0)s_1-As_1^*s_1 +s_2^*s_2A).
  $$

  If the range of $D\SF$ is not the whole tangent space then there
  exist $w\in TM$ and $B\in \Hom(\C^{p},\C^{n-p})$ with
  $<D\SF(v,A,b),(w,B)>=0$ for all $(v,A,b)$. Choose $v=0$, $A=0$ and
  $b(0)=0$ to get $\Tr(s_1^*\bd_w b(0) s_1)=0$. Varying $\bd b(0)$
  and using that $s_1\neq 0$ we get that $w=0$. Now choose $v=0$,
  $A=0$ to get $\Tr(B^* s_2^* b(0) s_1)=0$. Varying $b(0)$ among all
  hermitian matrices, we get $s_2B=0$. Now choose $v=0$ to get
  $\Tr(B^* As_1^*s_1)=0$ for all $A$. As $s_1\neq 0$ we get $B=0$,
  which completes the proof of the surjectivity of $D\SF$.

  Being $\SF$ a submersion, then for generic metric $g\in \SM_F$ the zero
  set of $\SF$ is regular. Therefore the set of critical points of
  $|\s|_g^2$ is a finite collection of points of $G$. The metric
  $g$ is not $\SC^\infty$ in principle but we may approximate it
  by a smooth metric without losing the above property.

  Also we want to check that the projection of the critical points
  to $M$ is in generic position. This is achieved if we check that
  at any $(y,g_0)$, for $v\in TM$ there exist $A\in
  \Hom(\C^p,\C^{n-p})$ and $b\in T\SM_F$ such that $D\SF(v,A,b)=0$.
  For this choose first $(A,b(0))$ such that $D\SF_2(v,A,b)=0$ and
  then for such $v$ and $A$ we can choose  $\bd b(0)$ such that
  $D\SF_1(v,A,b)=0$. Finally we take some $b \in
  T_{g_0}\SM$ with the obtained $b(0)$ and $\bd b(0)$.

  Now at a generic point of $M$ we have that $\p$ is injective.
  Therefore at the critical points of $|\s|^2_g$ the map $\p_r$ is
  injective, which is what we wanted to show.
\hfill $\Box$ \vspace{5mm}

{\bf Proof of Theorem \ref{thm:determ}.} By Proposition \ref{G2},
$(U^*\otimes F, \phi_r)$ is a $(k+r(e-r))$-positive pair. On the
other hand the rank of the bundle $U^*\otimes F$ is $(e-r)f$. By
Theorem \ref{main_thm} we get that $Z(\phi_r)$ satisfies  the
required isomorphisms with $i <k+r(e-r)-f(e-r)=k-(e-r)(f-r)$, and
an epimorphism for $i=k-(e-r)(f-r)$. \hfill $\Box$ \vspace{5mm}

\section{Application to Brill-Noether theory}\label{BN}

We briefly recall all the tools that we are going to use following
\cite{ACGH85}. We shall deal with the case of general rank since
many of the constructions are valid in this case, and later we
will particularize to the case of line bundles. Fix an algebraic
curve $C$ and denote by $\SM_d^r(C)$ the moduli space of stable
vector bundles of rank $r$ and degree $d$. We choose $d$ and $r$
coprime so that $\SM_d^r(C)$ is compact. Brill-Noether moduli
spaces, $W_d^{r,k}(C)$, are defined as
 $$
 W_d^{r,k}(C)= \{ A\in \SM_d^{r}(C)| ~ h^0(A)\geq k+1 \}.
 $$
It is also convenient to define the varieties
 $$
 G_d^{r,k}(C)= \{ (A,V) |A \in \SM_d^{r}(C), V \subset H^0(A),
 \dim V= k+1 \}.
 $$
The Brill-Noether loci for rank $r=1$ are the ones defined in the
introduction, with $\SM_d^{1}(C)=\Jac^d(C)$.

It is possible to choose a high degree divisor $\Gamma$ on $C$,
let us say of degree $m$, such that for each element $A\in
\SM_d^r$ the long exact sequence
 $$
 H^0(A) \to H^0(A(\G)) \stackrel{\phi}{\to} H^0(A(\G)/A)
  \to H^1(A) \to H^1(A(\G)) \to \cdots
 $$
satisfies that $H^1(A(\Gamma))=0$, since the bundles in
$\SM_d^r(C)$ form a bounded family. A bundle will be in
$W_d^{r,k}(C)$ whenever $H^0(A)$ is big enough. This is equivalent
to make $\rank\, \phi$ small enough. Denote $n=rm+d$. Take a
universal bundle $\SU$ over $\SM_d^r(C)\times C$ and consider the
sequence
 \begin{equation}
 \SU \to \SU(\Gamma) \to \SU(\Gamma)/\SU. \label{seq_uni}
 \end{equation}
If $\pi$ is the projection of $\SM_d^r\times C$ to $\SM_d^r$, we
obtain by pushing forward through $\pi$
 $$
 \pi_*(\SU(\Gamma)) \stackrel{\psi}{\longrightarrow}
 \pi_*(\SU(\Gamma)/\SU).
 $$
The determinantal varieties of $\psi$ define the sequence of
Brill-Noether loci $W_d^{r,k}(C)$. Moreover the spaces
$G_d^{r,k}(C)$ are the spaces $Z(\psi_s)$ in the notation of
Section \ref{determ}. So we need to study the curvature of the
vector bundles
 $$
  E=\pi_*(\SU(\Gamma)), \qquad F=\pi_*(\SU(\Gamma)/\SU)
 $$
to understand the topology of $W_d^{r,k}(C)$. More precisely,
$e=\rank E=d+mr+r(1-g)$ and $f=\rank F=mr$. Suppose for instance
that $\frac dr \leq g-1$. Then $e\leq f$ and the map $\psi$ is
generically injective. Put $s=e-(k+1)$. Then the locus
$Z(\psi_s)=G_d^{r,k}(C)\subset \Gr(e-s, E)$. The (expected)
codimension of $Z(\psi_s)$ is $(e-s)(f-s)=(k+1)(k+1-d+r(g-1))$.
The case $\frac dr > g-1$ is similar.

\begin{proposition} \label{Epositivo}
There is a (natural) metric on $E\to\SM_d^r(C)$ such that its
curvature is given by
 $$
 \la t, \Theta_{u,v}t\ra =-\ima \int_C \la G_{\End A}
 (\Lambda [\eta_u, \bar{\eta}_v]) t, t \ra \o
 + \int_C \la G_{A(\G)}(\eta_u t),\eta_v t\ra \o+\b(u,v)
  \int_C |t|^2\o,
 $$
at a point $A\in \SM_d^r(C)$, where $t\in
E_A=H^0(\SU(\G)|_{\{A\}\x C}) \isom H^0(A(\G))$, $u,v \in
T_A\SM_d^r(C)$. Here  $\o$ is the area form on $C$, $\b$ is a
$(1,1)$-form on $\SM_d^r(C)$, $\eta_u$ and $\eta_v$ are the
Kodaira-Spencer representatives, i.e.,\ the harmonic elements of
$\O^{0,1}(\End A)$ corresponding to $u,v$ under the natural
isomorphism $T_A\SM_d^r(C) \isom H^1(\End A)$, and $G_{\End A}$
and $G_{A(\G)}$ are the Green operators for the Laplacian
$\Delta_{\bbd}$ on $\O^0(\End A)$ and $\O^{0,1}(A(\G))$,
respectively.
\end{proposition}

{\bf Proof.} We write $M=\SM_d^r(C)$ for simplicity. Let us start
by computing the curvature of the universal bundle $\SU\to X=M\x
C$. Put a hermitian metric on $\SU$ so that it becomes a family of
hermitian-Einstein line bundles over $C$. Let $F$ stand for the
curvature of $\SU$ and decompose it as
$F=F_{CC}+F_{CM}+F_{MC}+F_{MM}$, according to the decomposition
$X=M\x C$. For instance, $F_{CM}\in \O_C^{1,0}\ox \O_M^{0,1}(\End
A)$, $F_{MM}\in \O^0_C\ox \O_M^{1,1}(\End A)$, etc. In particular
$F_{CC} =\l \o I$, where $\l$ is a constant.

Decompose the differential given by the natural connection on
$\SU$ as $d=d_C+d_M$. Clearly $\bbd_C F_{MC}=0$. Using the Bianchi
identity $dF=0$ and looking at the decomposition
 $$
 \O_X^{p,q}=\bigoplus_{i,j} \O_M^{i,j} \ox \O_C^{p-i,q-j},
 $$
we have that $\bd_C F_{MC}=-\bd_M F_{CC}=0$. Therefore
$F_{MC}|_{\{A\}\x C}$ is harmonic. Actually $F_{MC}$ gives a
complex linear map $T_AM \to \O^{0,1}_C(\End A)$, which represents
the Kodaira-Spencer map~\cite{TW}. Also
$F_{CM}=\overline{F}_{MC}$. Using again the Bianchi identity and
$\bd_C F_{MC}=0$ we have that
 $$
  \bd_C\bbd_C F_{MM}=- \bd_C \bbd_M F_{MC}= -\bbd_M \bd_C F_{MC}
  -[F_{CM}, F_{MC}]=- [F_{CM}, F_{MC}]=-[F_{MC},F_{CM}],
 $$
where this is a combination of the wedge product of forms and the
Lie bracket of the endomorphisms of $A$. By the K\"{a}hler
identities~\cite{We73} we have that $\Delta_{\bbd} F_{MM}= -\ima
\Lambda [F_{MC}, F_{CM}]$, so
 $$
 F_{MM} =-\ima G_{\End A} (\Lambda [F_{MC},F_{CM}])+\pi^*\b,
 $$
for a purely imaginary $(1,1)$-form $\b$ on $M$.

Now we compute the curvature of the push-forward $\pi_*
(\SU(\G))$. Since $\SU(\G)\to \SM_n^r(C)\x C$ is also a universal
bundle, where $n=rm+d$, we may suppose that $\G$ is the zero
divisor for this computation. The natural hermitian $L^2$-metric
$h_E$ on $E=\pi_*\SU$ is given by $h_E(t_1, t_2)= \int_C h(t_1,
t_2) \omega$, for $t_1,t_2\in E_A=H^0(A)$. Now the proof of
\cite[theorem 1]{TW} states that for the $L^2$ metric on $\pi_*
\SU$, we have for $u,v\in T_A M$ and $t \in H^0(A)$,
 \begin{eqnarray*}
 \la t,\Theta_{u,v}t \ra &= &\int_C \la F_{MM}(u,v))t,t\ra \o +\int_C
 \la G_{A}(\eta_u t),\eta_v t \ra\o = \\
 &=&  -\ima \int_C \la G_{\End A}(\Lambda [\eta_u, \bar{\eta}_v])
 t, t\ra\o
 +\int_C \la G_{A}(\eta_u t),\eta_v t \ra\o +\b(u,v)\int_C |t|^2\o,
 \end{eqnarray*}
since $\eta_u=F_{CM}(u)$. There is a different sign in our formula
to that of \cite[theorem 1]{TW} due to our convention on the
hermitian metric to be complex linear in the second variable.
\hfill $\Box$ \vspace{5mm}

In order to get more specific information on the bundles $E$ and
$F$ we need to restrict from now on to the case of line bundles,
i.e.,\ $r=1$. In this case $\SM^1_d(C)=\Jac^d(C)$ and the
Brill-Noether loci are those $W_d^k(C)$ described in the
introduction. We have the following

\begin{corollary}\label{rank1}
 Suppose that $r=1$. Then the
 bundle $E^*\otimes F$ on $\Jac^d(C)$ is positive.
\end{corollary}

{\bf Proof.}
 In this case the universal bundle is
the Poincar\'{e} bundle $\SL\to \Jac^d(C)\x C$. The bundle
 $$
  F= \pi_*(\SL(\Gamma)/\SL)=\bigoplus_{i=1}^m
  \SL|_{x_i\x \Jac^d(C)},
 $$
where $\G=x_1+\cdots+x_m$. So to check that $E^*\otimes F$ is
positive we only need to check that $E^*\otimes \SL|_{x_i\x
\Jac^d(C)}$ is so. This is equivalent to trivialize $\SL$ so that
$\SL|_{x_i\x \Jac^d(C)}\isom \SO$ and then to check that $E^*$ is
positive.

In the case of line bundles, the formula for the curvature of $E$
in Proposition \ref{Epositivo} reduces because the Lie bracket is
zero. Putting $v=\ima u$ we get for the curvature of $E$, for
$t\in E_A$ with $||t||=1$,
 $$
 \la t,(\Theta_E)_{u,\ima u} t\ra = \ima \int_C \la G_{A(\G)}(\eta_u t),
 \eta_u t \ra \o + \b(u,\ima u).
 $$

Moreover the curvature of $\SL$ in the Jacobian direction is then
$F_{MM}=\pi^*\b$. Restricting to $x_i\x \Jac^d(C)$ we get $\b$ is
an exact form. We may then change the hermitian metric in $E\to
\Jac^d(C)$ in order to eliminate this term from the formula of the
curvature. So we get $\la t,(\Theta_{E^*})_{u,\ima u} t\ra = -\ima
\int_C \la G_{A(\G)}(\eta_u t), \eta_u t \ra \o$. Now
$\Delta_{\bbd}$ is a positive operator on $\O^{0,1}(A(\G))$ since
it is self-adjoint semi-positive with kernel $H^1(A(\G))=0$, by
choosing the degree of $\G$ very large. Therefore $G_{A(\G)}$ is
also positive and the result follows. \hfill $\Box$ \vspace{5mm}

{\bf Proof of Theorem \ref{tue}.} By \cite[page 214]{ACGH85}, for
a general curve $C$ of genus $g$, $G_d^k(C)$ is reduced and of
pure dimension $\rho= g - (k+1)(g-d+k)$. Moreover $G_d^k(C)$ is
smooth.

Applying theorem \ref{thm:determ} using corollary \ref{rank1} we
get that the natural inclusion of $G_d^k(C)=Z(\psi_s)$ in
$G=\Gr(e-s, E)$ induces isomorphisms in the homology groups of
order less than or equal to $\rho-1$ and epimorphisms for order
$\rho$. The same holds for the homotopy groups.

Now $G$ is the total space of a fibration $\Gr(e-s,e)\to G\to
\Jac^d(C)=(\SSS^1)^{2g}$. Therefore $\pi_1(G)=\Z^{2g}$ and
$\pi_i(G)=\pi_i(\Gr(e-s,e))$ for $i>1$, since the grassmannian is
simply connected. Now $e-s=k+1$ and $e$ is very large (taking the
degree of $\G$ very large), so  $\Gr(e-s,e)=\Gr(k+1,e) \inc
\Gr(k+1,\infty) \cong BU(k+1)$ induces isomorphisms in homotopy
and homology groups up to some arbitrarily large $i_0$. This
proves the statements of the theorem regarding $G^k_d(C)$.

It only remains to prove that $\pi_1(G^k_d(C))=\pi_1(W^k_d(C))$.
Using $\pi_1(G)\isom \pi_1(G^k_d(C)) \to \pi_1(W^k_d(C)) \to
\pi_1(\Jac^d(C))$, we get injectivity on the fundamental groups.
We claim that the surjectivity follows from the fact that
$W_d^{k+i}(C)$, $i\geq 0$, form a Whitney stratification of
$W_d^k(C)$ and that the fibers of $\pi^{-1}(W_d^{k+i}(C)) \cap
G_d^k(C) \to W_d^{k+i}(C)$ are connected. In fact, given a
representative $\gamma:S^1 \to W^k_d(C)$ of a class in the
fundamental group, we can always assure that it intersects the
singular strata of $W^k_d(C)$ in isolated points $p_1, \ldots,
p_k$. Outside these points we can lift $\gamma$ to $\hat{\gamma}$
in $G^k_d(C)$. As the preimage of $p_i$ is connected,
$\hat{\gamma}$ can be extended to a loop, also denoted
$\hat{\gamma}$, mapping to $\gamma$. This completes the proof.
\hfill $\Box$ \vspace{5mm}

It is also interesting to study the higher homotopy groups of the
Brill-Noether varieties $W_d^k(C)$ in the case they are not
smooth, i.e.,\ when $W_d^{k+1}(C)\neq \emptyset$.

{\bf Proof of Corollary \ref{tue2}.} The idea is that
$W^{k+l}_d\neq\emptyset$ implies that there is a $\Gr(k+1,k+l+1)
\inc G_d^k$ mapping to a point in $W_d^k$. Since
$\pi_i(\Gr(k+1,k+l+1)) \stackrel{\simeq}{\to} \pi_i(G_d^k)\isom
\pi_i(\Gr(k+1,e))$ for $i\leq 2l$, we have that $\pi_i(G^k_d)\to
\pi_i(W^k_d)$ is the zero map. Instead of showing the surjectivity
of such map, we will take a shortcut. By \cite[theorem
2.2]{Debarre} the inclusion $W_d^k\inc \Jac^d$ induces isomorphism
on homology up to degree $2l$. By theorem \ref{tue} it also
induces isomorphism on the fundamental group. This implies that it
induces isomorphism on the homotopy groups up to degree $2l$.
\hfill $\Box$ 

\end{document}